\documentclass[10pt]{article}
 \font\smallit=cmti10

\usepackage{amssymb,latexsym,amsmath,epsfig,amsthm} 
\usepackage{enumitem}

\makeatletter

\renewcommand{\@seccntformat}[1]{\csname the#1\endcsname. }

\makeatother

 \newtheorem{theorem}{Theorem}[section]
 \newtheorem{lemma}[theorem]{Lemma}
 
 \newtheorem{proposition}[theorem]{Proposition}
 \newtheorem{corollary}[theorem]{Corollary}
 \newtheorem{definition}[theorem]{Definition}
 
 \newtheorem{example}[theorem]{Example}
 \newtheorem{main theorem}{Main Theorem}

\begin{document}
\begin{center}
 {\bf Norm Orthogonal Bases and Invariants of $p$-adic Lattices}
 \vskip 30pt

 {\bf Chi Zhang, Yingpu Deng and Zhaonan Wang}\\
 {\smallit State Key Laboratory of Mathematical Sciences, Academy of Mathematics and Systems Science, Chinese Academy of Sciences, Beijing 100190, People's Republic of China}\\
 {and}\\
 {\smallit School of Mathematical Sciences, University of Chinese Academy of Sciences, Beijing 100049, People's Republic of China}\\

 \vskip 10pt

 {\tt zhangchi171@mails.ucas.ac.cn, dengyp@amss.ac.cn, znwang@amss.ac.cn}\\

 \end{center}
\vskip 30pt

\centerline{\bf Abstract} In 2018, the Longest Vector Problem (LVP) and the Closest Vector Problem (CVP) in $p$-adic lattices were introduced. These problems are closely linked to the orthogonalization process. In this paper, we first prove that every $p$-adic lattice has an orthogonal basis respect to any given norm, whereas lattices in Euclidean spaces lack such bases in general. It is an improvement on Weil's result. Then, we prove that the sorted norm sequence of orthogonal basis of a $p$-adic lattice is unique and give definitions to the successive maxima and the escape distance, as the $p$-adic analogues of the successive minima and the covering radius in Euclidean lattices. Finally, we present deterministic polynomial time algorithms designed for the orthogonalization process, addressing both the LVP and the CVP with the help of an orthogonal basis of the whole vector space.\par


\noindent

\pagestyle{myheadings}

 \thispagestyle{empty}
 \baselineskip=12.875pt
 \vskip 20pt

\section{Introduction}

Lattices were investigated since the late $18$th century. There are two famous computational problems in Euclidean lattices: the Shortest Vector Problem (SVP) and the Closest Vector Problem (CVP). Much research has been done to study these problems. Van Emde Boas \cite{ref-6} proved that the SVP is NP-hard in the $l_{\infty}$ norm and the CVP is NP-hard in the $l_2$ norm by reducing them to the weak partition problem. Ajtai \cite{ref-7} proved that the SVP is NP-hard under randomized reductions in the $l_2$ norm. These hard problems are useful in constructing cryptographic primitives. The GGH scheme \cite{ref-8}, the NTRU scheme \cite{ref-9} and the LWE scheme \cite{ref-10} are based on these problems.\par
The $p$-adic numbers $\mathbb{Q}_p$ were invented by Hensel in the late 19th century. The concept of a local field is an abstraction of the field $\mathbb{Q}_p$. There is much research on the $p$-adic fields. It is well-known that a $p$-adic field has only finitely many extensions of a given finite degree. Krasner \cite{ref-12} gave formulas for the number of extensions of a given degree and discriminant. Following his work, Pauli et al. \cite{ref-13} presented an algorithm for the computation of generating polynomials for all extensions of a given degree and discriminant. Based on the Round Four maximal order algorithm \cite{ref-14}, Ford et al. \cite{ref-15} proposed an algorithm to factor polynomials over local fields. Later, Gu\`{a}rdia et al. \cite{ref-16} used the ``single-factor lifting" routine to provide a fast  polynomial factorization algorithm over local fields.\par
Lattices can also be defined in local fields such as $p$-adic fields, see \cite{ref-5}. However, there are not much research on properties and applications of the $p$-adic lattices. Motivated by lattice-based cryptosystems, one of the most promising post-quantum cryptosystems, Deng et al. \cite{ref-1} introduced two new computational problems in $p$-adic lattices of local fields, the Longest Vector Problem (LVP) and the Closest Vector Problem (CVP), which are the $p$-adic analogues of the Shortest Vector Problem and the Closest Vector Problem in lattices of Euclidean spaces. They considered these new problems to be challenging and potentially applicable for constructing public-key cryptosystems. Moreover, as the $p$-adic analogues of the lattices in Euclidean spaces, it is reasonable to expect these problems to be quantum-resistant. These new problems might contribute to the post-quantum assumptions.\par
Therefore, properties of the $p$-adic lattices are worth studying. Interestingly, $p$-adic lattices possess some totally different properties comparing with lattices in Euclidean spaces. For example, $\mathbb{Q}_p^n$ can be viewed as a field for every integer $n\ge1$, while the famous Frobenius Theorem asserts that $\mathbb{R}^n$ can be viewed as a field only when $n=1,2,4$. In section \ref{se-1}, we prove the following main theorem.
\begin{main theorem}\label{th-1.5}
Let $V$ be a vector space over $\mathbb{Q}_p$ of finite dimension, and let $N$ be a norm on $V$. Let $\alpha_{1},\dots,\alpha_{n}$ be $\mathbb{Q}_p$-linearly independent vectors of $V$. Let $\mathcal{L}=\mathcal{L}(\alpha_{1},\dots,\alpha_{n})$ be a $p$-adic lattice of rank $n$ in $V$. Then there is a deterministic algorithm to find out an $N$-orthogonal basis of the lattice $\mathcal{L}$ if we can compute efficiently the norm $N(v)$ of any vector $v\in V$.
\end{main theorem}
This theorem shows that every $p$-adic lattice has an orthogonal basis respect to any given norm, whereas lattices in Euclidean spaces lack such bases in general. Weil \cite{ref-5} proved that for every $p$-adic lattice, there exists a norm such that the $p$-adic lattice has an orthogonal basis respect to this norm. Our theorem improves Weil's result.\par
Additionally, we find an invariant of $p$-adic lattice which is called the successive maxima. This is the $p$-adic analogue of the successive minima in Euclidean lattices. We present another main theorem here and prove it in section \ref{se-2.1}.
\begin{main theorem}\label{th-2.1}
Let $V$ be a vector space over $\mathbb{Q}_p$ of finite dimension, and let $N$ be a norm on $V$. Let $\mathcal{L}$ be a $p$-adic lattice of rank $n$ in $V$. Suppose that $\alpha_{1},\dots,\alpha_{n}$ and $\beta_{1},\dots,\beta_{n}$ are two $N$-orthogonal bases of $\mathcal{L}$ such that $N(\alpha_{1})\ge\cdots\ge N(\alpha_{n})$ and $N(\beta_{1})\ge\cdots\ge N(\beta_{n})$. Then we have $N(\alpha_i)=N(\beta_i)$ for $1\le i\le n$.
\end{main theorem}

On the other hand, the LVP and the CVP are closely linked to the orthogonalization process. Deng et al. \cite{ref-2} proved that the LVP can be solved efficiently with  the help of an orthogonal basis of the $p$-adic lattice. Furthermore, if this orthogonal basis can be extended to an orthogonal basis of the whole vector space, then the CVP can be solved efficiently. We improve these algorithms in Section \ref{se-5}. The new algorithms only need an orthogonal basis of the whole vector space.\par
This paper is organized as follows. In Section \ref{se-0} we recall some definitions and properties. We prove our Main Theorem 1 in Section \ref{se-1}. Subsequently, we prove our Main Theorem 2 in Section \ref{se-2}. We also give the definitions to the successive maxima and the escape distance in Section \ref{se-2}. In Section \ref{se-3} we prove that any two orthogonal bases of a $p$-adic lattice can be obtained from each other by some operations. Then, we present a deterministic polynomial time algorithm to find orthogonal bases of $p$-adic lattices with the help of an orthogonal basis of the whole vector space in Section \ref{se-4}. Next, we introduce new deterministic polynomial time algorithms to solve the CVP and the LVP with the help of orthogonal bases in Section \ref{se-5}. Finally, we conclude that the orthogonalization process and the CVP are polynomially equivalent.

\section{Preliminaries}\label{se-0}
\subsection{Norm and Orthogonal Basis}

Let $p$ be a prime. Let $V$ be a vector space over $\mathbb{Q}_p$ of finite dimension. A norm $N$ on $V$ is a function
$$N:V\rightarrow\mathbb{R}$$
such that

\begin{enumerate}
\item $N(v)\ge0$ for any $v\in V$, and $N(v)=0$ if and only if $v=0$;
\item $N(xv)=\left|x\right|_{p}\cdot N(v)$ for any $x\in\mathbb{Q}_p$ and $v\in V$;
\item $N(v+w)\le\max{\left\{N(v),N(w)\right\}}$ for any $v,w\in V$.
\end{enumerate}

\noindent Here, $\left|x\right|_{p}$ is the $p$-adic absolute value for a $p$-adic number $x\in\mathbb{Q}_p$.\par

If $N$ is a norm on $V$, and if $N(v)\ne N(w)$ for $v,w\in V$, then we must have $N(v+w)=\max{\{N(v),N(w)\}}$. Weil (\cite{ref-5} page 26) proved the following proposition.

\begin{proposition}[\cite{ref-5}]
Let $V$ be a vector space over $\mathbb{Q}_p$ of finite dimension $n>0$, and let $N$ be a norm on $V$. Then there is a decomposition $V=V_1+\cdots+V_n$ of $V$ into a direct sum of subspaces $V_i$ of dimension $1$, such that
$$N\left(\sum_{i=1}^{n}{v_i}\right)=\max_{1\le i\le n}{N(v_i)}$$
for any $v_i\in V_i$, $i=1,\dots,n$.
\end{proposition}

Weil proved the above proposition for finite-dimensional vector spaces over a $p$-field (commutative or not). For simplicity, we only consider the case $\mathbb{Q}_p$. Thus, we can define the orthogonal basis.

\begin{definition}[$N$-orthogonal basis \cite{ref-5}]
Let $V$ be a vector space over $\mathbb{Q}_p$ of finite dimension $n>0$, and let $N$ be a norm on $V$. We call $\alpha_1,\dots,\alpha_n$ an $N$-orthogonal basis of $V$ over $\mathbb{Q}_p$ if $V$ can be decomposed into the direct sum of $n$ $1$-dimensional subspaces $V_i$'s $(1\le i\le n)$, such that
$$N\left(\sum_{i=1}^{n}{v_i}\right)=\max_{1\le i\le n}{N\left(v_i\right)}$$
for any $v_i\in V_i$, $i=1,\dots,n$, where $V_i$ is spanned by $\alpha_i$.
Two subspaces $U$, $W$ of $V$ are said to be $N$-orthogonal if the sum $U+W$ is a direct sum and it holds that $N\left(u+w\right)=\max\left\{N(u), N(w)\right\}$ for all $u\in U$, $w\in W$.
\end{definition}

\subsection{$p$-adic Lattice}

We first recall the definition of a $p$-adic lattice.

\begin{definition}[$p$-adic lattice \cite{ref-1}]
Let $V$ be a vector space over $\mathbb{Q}_p$ of finite dimension $n>0$, and let $N$ be a norm on $V$. Let $\alpha_1,\dots,\alpha_m$ $(1\le m\le n)$ be $\mathbb{Q}_p$-linearly independent vectors of $V$. A $p$-adic lattice in $V$ is the set
$$\mathcal{L}(\alpha_1,\dots,\alpha_m):=\left\{\sum^{m}_{i=1}{a_i\alpha_i}:a_i\in\mathbb{Z}_p,1\le i\le m\right\}$$
of all $\mathbb{Z}_p$-linear combinations of $\alpha_1,\dots,\alpha_m$. The sequence of vectors $\alpha_1,\dots,\alpha_m$ is called a basis of the lattice $\mathcal{L}(\alpha_1,\dots,\alpha_m)$. The integers $m$ and $n$ are called the rank and dimension of the lattice, respectively. When $n=m$, we say that the lattice is of full rank.
\end{definition}

$p$-adic lattices are compact subsets of $V$. The following proposition can be found in \cite{ref-4} page 72, prop.

\begin{proposition}[\cite{ref-4}]\label{pr-0.4}
Let $\Omega\subset V$ be a compact subset.
\begin{enumerate}[label=(\alph*)]
\item For every $a\in V\setminus\Omega$, the set of norms $\{N(x-a):x\in\Omega\}$ is finite.
\item For every $a\in\Omega$, the set of norms $\{N(x-a):x\in\Omega\setminus\{a\}\}$ is discrete in $\mathbb{R}_{>0}$.
\end{enumerate}
\end{proposition}

We can also define the orthogonal basis of a $p$-adic lattice.

\begin{definition}[$N$-orthogonal basis of a $p$-adic lattice \cite{ref-3}]
Let $V$ be a vector space over $\mathbb{Q}_p$ of finite dimension $n>0$, and let $N$ be a norm on $V$. If $\alpha_1,\dots,\alpha_m$ is an $N$-orthogonal basis of the vector space spanned by a $p$-adic lattice $\mathcal{L}=\sum_{i=1}^{m}{\mathbb{Z}_p\alpha_i}$, then we call $\alpha_1,\dots,\alpha_m$ an $N$-orthogonal basis of the lattice $\mathcal{L}$.
\end{definition}

\subsection{LVP and CVP}

Deng et al. \cite{ref-1} introduced two new computational problems in $p$-adic lattices. They are the Longest Vector Problem (LVP) and the Closest Vector Problem (CVP). We first review them briefly.

\begin{definition}[\cite{ref-1}]
Let $\mathcal{L}=\mathcal{L}(\alpha_1,\dots,\alpha_m)$ be a $p$-adic lattice in $V$. We define recursively a sequence of positive real numbers $\lambda_1,\lambda_2,\lambda_3,\dots$ as follows.
$$\lambda_1=\max_{1\le i\le m}{N(\alpha_i)},$$
$$\lambda_{j+1}=\max{\left\{N(v):v\in\mathcal{L},N(v)<\lambda_j\right\}} \mbox{ for } j\ge1.$$
\end{definition}

We have $\lambda_1>\lambda_2>\lambda_3>\dots$ and $\lim_{j\rightarrow\infty}\lambda_j=0$. The Longest Vector Problem is defined as follows.

\begin{definition}[\cite{ref-1}]
Given a $p$-adic lattice $\mathcal{L}=\mathcal{L}(\alpha_1,\dots,\alpha_m)$ in $V$, the Longest Vector Problem is to find a lattice vector $v\in\mathcal{L}$ such that $N(v)=\lambda_2$.
\end{definition}

The Closest Vector Problem is defined as follows.

\begin{definition}[\cite{ref-1}]
Let $\mathcal{L}=\mathcal{L}(\alpha_1,\dots,\alpha_m)$ be a $p$-adic lattice in $V$ and let $t\in V$ be a target vector. The Closest Vector Problem is to find a lattice vector $v\in\mathcal{L}$ such that $N(t-v)=\min_{w\in\mathcal{L}}{N(t-w)}$.
\end{definition}

Thanks to Proposition \ref{pr-0.4}, these definitions are well-defined. Deng et al. \cite{ref-1} provided deterministic exponential time algorithms to solve the LVP and the CVP. Additionally, Deng et al. \cite{ref-2} presented deterministic polynomial time algorithms for solving the LVP and the CVP specifically with the help of orthogonal bases.

\section{Orthogonalization of $p$-adic Lattices}\label{se-1}

In this section, we prove that every $p$-adic lattice has an $N$-orthogonal basis. Firstly, we need some lemmas. Lemma \ref{le-1.1}, Corollary \ref{co-1.2} and Lemma \ref{le-1.3} can be found in \cite{ref-3}.

\begin{lemma}[\cite{ref-3}]\label{le-1.1}
Let $V$ be a vector space over $\mathbb{Q}_p$ of finite dimension, and let $N$ be a norm on $V$. Let $v,w\in V$. Then we have $N(v+w)=\max\{N(v),N(w)\}$ if and only if $N(v+w)\ge N(v)$.
\end{lemma}

\begin{corollary}[\cite{ref-3}]\label{co-1.2}
Let $V$ be a vector space over $\mathbb{Q}_p$ of finite dimension, and let $N$ be a norm on $V$. Let $\alpha_{1},\dots,\alpha_{n}$ $(n>1)$ be $\mathbb{Q}_p$-linearly independent vectors of $V$. Set $\mathcal{L}=\mathcal{L}(\alpha_{2},\dots,\alpha_{n})$. Then we have:
$$N(\alpha_{1}+w)=\max\{N(\alpha_{1}),N(w)\} \mbox{ for all } w\in \mathcal{L}$$
if and only if
$$N(\alpha_{1})=\min\{N(\alpha_{1}+w) : w\in\mathcal{L}\}.$$
\end{corollary}

\begin{lemma}[\cite{ref-3}]\label{le-1.3}
Let $V$ be a vector space over $\mathbb{Q}_p$ of finite dimension $n>1$, and let $N$ be a norm on $V$. Let $\alpha_{1},\dots,\alpha_{n}$ be $\mathbb{Q}_p$-linearly independent vectors of $V$. Then $\alpha_{1},\dots,\alpha_{n}$ is an $N$-orthogonal basis of $V$ if and only if it holds that
$$N\left(\sum^{n}_{i=1}a_{i}\alpha_{i}\right)=\max_{1\le i\le n}N(a_{i}\alpha_{i}),$$
where one of the $a_{1},\dots,a_{n}$ is $1$ and the others are in $\mathbb{Z}_p$.
\end{lemma}

The next lemma shows that we can find an $N$-orthogonal basis of a $p$-adic lattice $\mathcal{L}$ by solving a CVP-instance if $\mathcal{L}$ has some special structure.

\begin{lemma}\label{le-1.4}
Let $V$ be a vector space over $\mathbb{Q}_p$ of finite dimension, and let $N$ be a norm on $V$. Let $\alpha_{1},\dots,\alpha_{n}$ $(n>1)$ be $\mathbb{Q}_p$-linearly independent vectors of $V$. Let $\mathcal{L}=\mathcal{L}(\alpha_{1},\dots,\alpha_{n})$ be a $p$-adic lattice of rank $n$ in $V$. Assume that $\alpha_{1},\dots,\alpha_{n-1}$ is an $N$-orthogonal basis of the lattice $\mathcal{L}(\alpha_{1},\dots,\alpha_{n-1})$ and $N(\alpha_{n})\le N(\alpha_{i})$ for $1\le i\le n-1$. Let $w_{0}\in\mathcal{L}(\alpha_{1},\dots,\alpha_{n-1})$ be such that
$$N(\alpha_{n}+w_{0})=\min\{N(\alpha_{n}+w):w\in\mathcal{L}(\alpha_{1},\dots,\alpha_{n-1})\}.$$
Let $\alpha^{\prime}_{n}=\alpha_{n}+w_{0}$. Then $\alpha_{1},\dots,\alpha_{n-1},\alpha^{\prime}_{n}$ is an $N$-orthogonal basis of $\mathcal{L}$.
\begin{proof}
Since $w_{0}\in\mathcal{L}(\alpha_{1},\dots,\alpha_{n-1})$ and $\alpha^{\prime}_{n}=\alpha_{n}+w_{0}$, we have
$$\mathcal{L}=\mathcal{L}(\alpha_{1},\dots,\alpha_{n-1},\alpha^{\prime}_{n}).$$
Next, we prove that $\alpha_{1},\dots,\alpha_{n-1},\alpha^{\prime}_{n}$ is an $N$-orthogonal basis of $\mathcal{L}$.\par
By Lemma \ref{le-1.3}, we only need to consider vectors $\sum^{n-1}_{i=1}a_{i}\alpha_{i}+a_{n}\alpha^{\prime}_{n}$ such that one of the $a_{1},\dots,a_{n}$ is $1$ and the others are in $\mathbb{Z}_p$. If $a_{n}\in\mathbb{Z}_p\setminus p\mathbb{Z}_p$, then $\left|a_{n}\right|_{p}=1$, hence it is invertible. By Corollary \ref{co-1.2}, we have
\begin{equation*}
\begin{split}
N\left(\sum^{n-1}_{i=1}a_{i}\alpha_{i}+a_{n}\alpha^{\prime}_{n}\right)&=N\left(\sum^{n-1}_{i=1}\frac{a_{i}}{a_{n}}\alpha_{i}+\alpha^{\prime}_{n}\right)\\
&=\max\left\{N\left(\sum^{n-1}_{i=1}\frac{a_{i}}{a_{n}}\alpha_{i}\right),N(\alpha^{\prime}_{n})\right\}\\
&=\max\left\{N\left(\sum^{n-1}_{i=1}a_{i}\alpha_{i}\right),N(a_{n}\alpha^{\prime}_{n})\right\}\\
&=\max\{N\left(a_{1}\alpha_{1}),\dots,N(a_{n-1}\alpha_{n-1}\right),N(a_{n}\alpha^{\prime}_{n}))\}.
\end{split}
\end{equation*}
On the other hand, if $a_{n}\in p\mathbb{Z}_p$, then $\left|a_{n}\right|_{p}<1$ and there is a $1$ among $a_{1},\dots,a_{n-1}$. We may assume that $a_{1}=1$. Since
$$N\left(\sum^{n-1}_{i=1}a_{i}\alpha_{i}\right)=\max_{1\le i\le n-1}N(a_{i}\alpha_{i})\ge N(\alpha_{1})>N(a_{n}\alpha_{n})\ge N(a_{n}\alpha^{\prime}_{n}),$$
we have 
\begin{equation*}
\begin{split}
N\left(\sum^{n-1}_{i=1}a_{i}\alpha_{i}+a_{n}\alpha^{\prime}_{n}\right)&=N\left(\sum^{n-1}_{i=1}a_{i}\alpha_{i}\right)\\
&=\max\{N(a_{1}\alpha_{1}),\dots,N(a_{n-1}\alpha_{n-1}),N(a_{n}\alpha^{\prime}_{n})\}.
\end{split}
\end{equation*}
Hence $\alpha_{1},\dots,\alpha_{n-1},\alpha^{\prime}_{n}$ is an $N$-orthogonal basis of $\mathcal{L}$.
\end{proof}
\end{lemma}

Now we can prove the Main Theorem \ref{th-1.5}. This result gives an affirmative answer to the question in \cite{ref-3}, asking whether $p$-adic lattices with rank greater than or equal to three have orthogonal bases.

\begin{proof}
Define $CVP(t,\mathcal{L})=w_0$ such that $$N(t-w_0)=\min\{N(t-w):w\in\mathcal{L}\}.$$ 
For instance, Theorem 4.4 in \cite{ref-1} provides such an algorithm. We will also give a CVP algorithm later in Section \ref{se-5}. The algorithm for performing the orthogonalization process is presented as follows. \par\vskip 10pt

{\bf Algorithm} (orthogonalization).\par
{\bf Input:} a basis $\alpha_{1},\dots,\alpha_{n}$ of a $p$-adic lattice $\mathcal{L}$.\par
{\bf Output:} an $N$-orthogonal basis of $\mathcal{L}$.
\begin{enumerate}
\item for $i=1$ to $n$ do:
\item\quad rearrange $\alpha_{i},\dots,\alpha_{n}$ such that  $N(\alpha_{i})=\max_{i\le k\le n}{N(\alpha_{k})}$,
\item\quad for $j=i+1$ to $n$ do:
\item\quad\quad $\alpha_{j}\leftarrow\alpha_{j}-CVP(\alpha_{j},\mathcal{L}(\alpha_{1},\dots,\alpha_{i})).$
\end{enumerate}\par
Return $(\alpha_{1},\dots,\alpha_{n})$.\par\vskip 10pt

This algorithm calls the CVP algorithm $O(n^2)$ times. We may assume that the order of $\alpha_{i},\dots,\alpha_{n}$ never changes in the step $2$. This can be done by giving the input vectors a proper order initially. For example, let the input be $\alpha_1,\alpha_2,\alpha_3$. Suppose that $N(\alpha_2)=\max{\{N(\alpha_1),N(\alpha_2),N(\alpha_3)\}}$, $CVP(\alpha_1,\mathcal{L}(\alpha_2))=w_1$, $CVP(\alpha_3,\mathcal{L}(\alpha_2))=w_3$ and $N(\alpha_3-w_3)\ge N(\alpha_1-w_1)$, then a proper order is $\alpha_2,\alpha_3,\alpha_1$.\par
Moreover, based on this assumption, the vector $\alpha_{i}$ will be fixed after the $(i-1)$th iteration of the outer loop (where the $0$th iteration of the outer loop means the beginning of the algorithm). In order to distinguish from the input vectors, let us denote this vector as $\beta_{i}$, then the input lattice of the CVP algorithm in the step $4$ will be $\mathcal{L}(\beta_{1},\dots,\beta_{i})$ and the output will be $\beta_{1},\dots,\beta_{n}$. \par
First, we prove that $N(\beta_{1})\ge\cdots\ge N(\beta_{n})$. For any $1\le i\le n-1$, suppose that we have obtained $\beta_{1},\dots,\beta_{i},\alpha^{\prime}_{i+1},\dots,\alpha^{\prime}_{n}$ after the $(i-1)$th iteration of the outer loop. According to the above assumptions, we have  $N(\beta_{i})\ge N(\alpha^{\prime}_{i+1})$. Since
$$\beta_{i+1}=\alpha^{\prime}_{i+1}-CVP(\alpha^{\prime}_{i+1},\mathcal{L}(\beta_{1},\dots,\beta_{i})),$$
we have $N(\beta_{i})\ge N(\alpha^{\prime}_{i+1})\ge N(\beta_{i+1})$. Hence $N(\beta_{1})\ge\cdots\ge N(\beta_{n})$.\par
Next, let us prove by induction that $\beta_{1},\dots,\beta_{i}$ is an $N$-orthogonal basis of $\mathcal{L}(\alpha_{1},\dots,\alpha_{i})$ $(1\le i\le n)$. When $i=1$, according to our assumptions, we have $\beta_{1}=\alpha_{1}$. It is clear that $\beta_{1}$ is an $N$-orthogonal basis of $\mathcal{L}(\beta_{1})=\mathcal{L}(\alpha_{1})$. Suppose that the conclusion holds for $i=k$, then $\beta_{1},\dots,\beta_{k}$ is an $N$-orthogonal basis of $\mathcal{L}(\alpha_{1},\dots,\alpha_{k})$ by induction hypothesis. Suppose that we have obtained $\beta_{1},\dots,\beta_{k},\alpha^{\prime\prime}_{k+1},\dots,\alpha^{\prime\prime}_{n}$ after the $(k-1)$th iteration of the outer loop. Since
$$N(\beta_{1})\ge\cdots\ge N(\beta_{k})\ge N(\alpha^{\prime\prime}_{k+1})$$
and
$$\beta_{k+1}=\alpha^{\prime\prime}_{k+1}-CVP(\alpha^{\prime\prime}_{k+1},\mathcal{L}(\beta_{1},\dots,\beta_{k})),$$
by Lemma \ref{le-1.4}, we conclude that $\beta_{1},\dots,\beta_{k+1}$ is an $N$-orthogonal basis of
$$\mathcal{L}(\beta_{1},\dots,\beta_{k},\alpha^{\prime\prime}_{k+1})=\mathcal{L}(\beta_{1},\dots,\beta_{k},\alpha_{k+1})=\mathcal{L}(\alpha_{1},\dots,\alpha_{k},\alpha_{k+1}),$$
where the first equality holds because $\alpha^{\prime\prime}_{k+1}=\alpha_{k+1}-w$ for some $w\in\mathcal{L}(\beta_{1},\dots,\beta_{k})$. Hence $\beta_{1},\dots,\beta_{n}$ is an $N$-orthogonal basis of $\mathcal{L}(\alpha_{1},\dots,\alpha_{n})$.
\end{proof}

Here is a toy example for illustrating the orthogonalization process in the proof of Main Theorem \ref{th-1.5}.

\begin{example}\label{eg-1.6}
Let $V=\mathbb{Q}_2(\zeta)$ where $\zeta$ is a primitive $5$th root of unity. Let $N$ be the unique absolute value extended by the $2$-adic absolute value of $\mathbb{Q}_p$ (see \cite{ref-4}). Since
$$N(a_1+a_2\zeta+a_3\zeta^2+a_4\zeta^3)=\max\{N(a_1),N(a_2\zeta),N(a_3\zeta^2),N(a_4\zeta^3)\},$$
where one of the $a_{1},\dots,a_{4}$ is $1$ and the others are in $\mathbb{Z}_p$. We conclude that $1,\zeta,\zeta^2,\zeta^3$ is an $N$-orthogonal basis of $V$ over $\mathbb{Q}_p$ by Lemma \ref{le-1.3}.\par
Let
$$\mathcal{L}=\mathcal{L}(1,1+2\zeta,2+8\zeta+16\zeta^2+16\zeta^3).$$
Since $N(1)=N(1+2\zeta)=1$ and $N(2+8\zeta+16\zeta^2+16\zeta^3)=\frac{1}{2}$, by solving the CVP-instances with the lattice $\mathcal{L}(1)$ and the target vectors $1+2\zeta$ and $2+8\zeta+16\zeta^2+16\zeta^3$, we obtain closest vectors $1$ and $2$, respectively. Notice that closest vector is not unique in general. Here the lattice vector $3$ in the lattice $\mathcal{L}(1)$ is also a closest vector of the target vector $1+2\zeta$.\par
Next, since $N(2\zeta)=\frac{1}{2}$ and $N(8\zeta+16\zeta^2+16\zeta^3)=\frac{1}{8}$, by solving the CVP-instance with the lattice $\mathcal{L}(1,2\zeta)$ and the target vector $8\zeta+16\zeta^2+16\zeta^3$, we obtain a closest vector $8\zeta$. Hence $1,2\zeta,16\zeta^2+16\zeta^3$ is an $N$-orthogonal basis of $\mathcal{L}$.
\end{example}

\section{Successive Maxima}\label{se-2}

In the realm of Euclidean lattices, the notion of successive minima is an important concept, which denotes the lengths of the shortest yet linearly independent vectors. We aim to explore the $p$-adic analogue of this concept. However, a straightforward substitution of ``minima'' with ``maxima'' results in a trivial definition, which does not capture the essence of the concept in the $p$-adic context. Consider a vector space $V$ over $\mathbb{Q}_p$ and a $p$-adic lattice $\mathcal{L}=\mathcal{L}(\alpha_1,\dots,\alpha_n)$ of rank $n$ in $V$ such that $\alpha_1,\dots,\alpha_n$ is an $N$-orthogonal basis of $\mathcal{L}$ and $N(\alpha_1)>\cdots>N(\alpha_n)$. Then the longest vector can be chosen as $\alpha_1$.\par
However, we can choose $\alpha_1+\alpha_2,\dots,\alpha_1+\alpha_n$ which are linearly independent over $\mathbb{Q}_p$ and all possess a common length of $N(\alpha_1)$. In order to avoid this trivial situation, we may require the vectors to be $N$-orthogonal rather than just linearly independent. This introduces the problem of determining the uniqueness of the sorted norm sequence of $N$-orthogonal bases of a $p$-adic lattice.\par
In this section, we firstly establish the uniqueness of the sorted norm sequence for $N$-orthogonal bases of a $p$-adic lattice. Subsequently, we provide definitions for the successive maxima and escape distance. Finally, we prove several properties associated with them.

\subsection{Proof of the Uniqueness}\label{se-2.1}

Let us begin with a simple observation.

\begin{proposition}\label{pr-2.0}
Let $V$ be a vector space over $\mathbb{Q}_p$ of finite dimension, and let $N$ be a norm on $V$. Let $\mathcal{L}$ be a $p$-adic lattice of rank $2$ in $V$. Suppose that $\alpha_{1},\alpha_{2}$ and $\beta_{1},\beta_{2}$ are two $N$-orthogonal bases of $\mathcal{L}$ such that $N(\alpha_{1})\ge N(\alpha_{2})$ and $N(\beta_{1})\ge N(\beta_{2})$. Then we have $N(\alpha_1)=N(\beta_1)$ and $N(\alpha_2)=N(\beta_2)$.
\begin{proof}
If $N(\alpha_1)\ne N(\beta_1)$, then we may assume that $N(\alpha_1)<N(\beta_1)$. Write
$$\beta_1=a_{11}\alpha_1+a_{12}\alpha_2,$$
where $a_{11},a_{12}\in\mathbb{Z}_p$. Since $\alpha_{1},\alpha_{2}$ is an $N$-orthogonal basis of $\mathcal{L}$, we have
$$N(\beta_1)=\max{\{N(a_{11}\alpha_1),N(a_{12}\alpha_2)\}}\le N(\alpha_1).$$
This is a contradiction. Hence $N(\alpha_1)=N(\beta_1)$.\par
If $N(\alpha_2)\ne N(\beta_2)$, then we may assume that $N(\alpha_2)<N(\beta_2)$. Now we have
$$N(\alpha_1)=N(\beta_1)\ge N(\beta_2)>N(\alpha_2).$$
Write
$$\beta_1=a_{11}\alpha_1+a_{12}\alpha_2,$$
$$\beta_2=a_{21}\alpha_1+a_{22}\alpha_2,$$
where $a_{11},a_{12},a_{21},a_{22}\in\mathbb{Z}_p$. Since $\alpha_{1},\alpha_{2}$ is an $N$-orthogonal basis of $\mathcal{L}$ and $N(\beta_1)=N(\alpha_1)>N(\alpha_2)$, we have $a_{11}\not\in p\mathbb{Z}_p$. Therefore,
$$N(\beta_2-a_{21}a^{-1}_{11}\beta_1)=N\big((a_{22}-a_{12}a_{21}a^{-1}_{11})\alpha_2\big)\le N(\alpha_2).$$
On the other hand, since $\beta_{1},\beta_{2}$ is an $N$-orthogonal basis of $\mathcal{L}$, we have
$$N(\beta_2-a_{21}a^{-1}_{11}\beta_1)=\max{\left\{N(\beta_2),N(a_{21}a^{-1}_{11}\beta_1)\right\}}\ge N(\beta_2).$$
This is a contradiction. Hence $N(\alpha_2)=N(\beta_2)$.
\end{proof}
\end{proposition}

Now we prove our Main Theorem \ref{th-2.1}. We prove by induction that if $N(\alpha_1)=N(\beta_1),N(\alpha_2)=N(\beta_2),\dots,N(\alpha_i)=N(\beta_i)$, then $N(\alpha_{i+1})=N(\beta_{i+1})$. The crucial point is that we can transform one $N$-orthogonal basis of $\mathcal{L}$ to another by a sequence of elementary row operations.

\begin{proof}
If $N(\alpha_1)<N(\beta_1)$, then we can write
$$\beta_1=\sum^{n}_{j=1}{a_{1j}\alpha_j},$$
where $a_{1j}\in\mathbb{Z}_p$ for $1\le j\le n$. Since $\alpha_{1},\dots,\alpha_{n}$ is an $N$-orthogonal basis of $\mathcal{L}$, we have
$$N(\beta_1)=\max_{1\le j\le n}{N(a_{1j}\alpha_j)}\le N(\alpha_1).$$
This is a contradiction.\par
If $N(\alpha_j)=N(\beta_j)$ for $1\le j\le m<n$ and $N(\alpha_{m+1})<N(\beta_{m+1})$,  then we can write
$$\beta_i=\sum^{n}_{j=1}{a_{ij}\alpha_j},$$
where $a_{ij}\in\mathbb{Z}_p$ for $1\le i\le m+1$ and $1\le j\le n$. Suppose that for $i=s_1,s_2,\dots,s_k\le m$, we have $N(\alpha_i)>N(\alpha_{i+1})$ and $N(\alpha_i)=N(\alpha_{i+1})$ for the remaining $1\le i\le m$. Set $s_0=0$, then for $1\le t\le k$, we can conclude that $a_{ij}\in p\mathbb{Z}_p$ if $s_{t-1}+1\le i\le s_t$ and $1\le j\le s_{t-1}$. Otherwise, there exists an $a_{i^{\prime}j^{\prime}}\not\in p\mathbb{Z}_p$ whose indices satisfy the above condition. Hence,
\begin{equation*}
\begin{split}
N(\alpha_{i^{\prime}})=N(\beta_{i^{\prime}})&=N\left(\sum^{n}_{j=1}{a_{i^{\prime}j}\alpha_j}\right)\\
&=\max_{1\le j\le n}{N(a_{i^{\prime}j}\alpha_j)}\ge N(a_{i^{\prime}j^{\prime}}\alpha_{j^{\prime}})=N(\alpha_{j^{\prime}})>N(\alpha_{i^{\prime}}).
\end{split}
\end{equation*}
This is a contradiction.\par
Let $u_i$ be the $(s_t-s_{t-1})$-tuple $(\overline{a_{i(s_{t-1}+1)}},\overline{a_{i(s_{t-1}+2)}},\dots,\overline{a_{is_{t}}})$ for $s_{t-1}+1\le i\le s_{t}$, $1\le t\le k$, where $\overline{a_{ij}}$ denote the image of $a_{ij}$ in $\mathbb{Z}_p/p\mathbb{Z}_p$. We claim that $u_{s_{t-1}+1},u_{s_{t-1}+2},\dots,u_{s_t}$ are linearly independent over $\mathbb{Z}_p/p\mathbb{Z}_p$ for each $1\le t\le k$. We prove this claim by induction.\par
First, when $t=1$, if there exist coefficients $\overline{b_1},\overline{b_2},\dots,\overline{b_{s_1}}\in\mathbb{Z}_p/p\mathbb{Z}_p$, not all $\overline{0}$, such that
$$\overline{b_1}u_1+\overline{b_2}u_2+\cdots+\overline{b_{s_1}}u_{s_1}=\overline{0},$$
then the coefficients of $\alpha_j$ $(1\le j\le s_1)$ in $\sum^{s_1}_{i=1}{b_i\beta_i}$ all belong to $p\mathbb{Z}_p$, where $b_i$ is an arbitrary pullback of $\overline{b_i}$. We can write
$$\sum^{s_1}_{i=1}{b_i\beta_i}=\sum^{n}_{i=1}{c_i\alpha_i}$$
for some coefficients $c_1,\dots,c_n\in\mathbb{Z}_p$. Hence
$$N\left(\sum^{s_1}_{i=1}{b_i\beta_i}\right)=N\left(\sum^{n}_{i=1}{c_i\alpha_i}\right)=\max_{1\le i\le n}{N(c_i\alpha_i)}<N(\alpha_1)=N(\beta_1).$$
On the other hand, since $\beta_{1},\dots,\beta_{n}$ is also an $N$-orthogonal basis of $\mathcal{L}$ and not all $b_1,\dots,b_{s_1}$ belong to $p\mathbb{Z}_p$, we have
$$N\left(\sum^{s_1}_{i=1}{b_i\beta_i}\right)=\max_{1\le i\le s_1}{N(b_i\beta_i)}=N(\beta_1).$$
This is a contradiction.\par
Moreover, we conclude that the matrix
$$\begin{pmatrix}
\overline{a_{11}} & \cdots & \overline{a_{1s_1}} \\
\vdots & \ddots & \vdots \\
\overline{a_{s_11}} & \cdots & \overline{a_{s_1s_1}}
\end{pmatrix}$$
is invertible over $\mathbb{Z}_p/p\mathbb{Z}_p$. Consequently, it can be transformed into an identity matrix by a sequence of elementary row operations. The preimage of this identity matrix is of the form
$$\begin{pmatrix}
a^{\prime}_{11} & \cdots & p\mathbb{Z}_p \\
\vdots & \ddots & \vdots \\
p\mathbb{Z}_p & \cdots & a^{\prime}_{s_1s_1}
\end{pmatrix},$$
where the diagonal elements are not in $p\mathbb{Z}_p$, while all other entries belong to $p\mathbb{Z}_p$. Hence we can multiply the first row by $(a^{\prime}_{11})^{-1}$ and use it to eliminate the other entries of the first column. It can be observed that, during the execution of these elementary row transformations, we consistently multiply the first row by multiples of $p$ and add the result to other rows. Consequently, even after these transformations, the diagonal elements still do not belong to $p\mathbb{Z}_p$. Hence, we can repeat the same procedure for the second row and proceed iteratively. Finally, we get an identity matrix, which means that we obtain a new basis $\beta^{\prime}_1,\dots,\beta^{\prime}_{s_1}$ of $\mathcal{L}(\beta_1,\dots,\beta_{s_1})$. The first $s_1$ coefficients of $\beta^{\prime}_1,\dots,\beta^{\prime}_{s_1}$ make up an identity matrix.\par
Now we can continue our induction. Suppose that $u_{s_{t-1}+1},u_{s_{t-1}+2},\dots,u_{s_t}$ are linearly independent over $\mathbb{Z}_p/p\mathbb{Z}_p$ for each $1\le t\le l-1$. According to the induction hypothesis and the previous conclusion, there is a basis of $\mathcal{L}(b_1,\dots,b_{s_{l-1}})$ of the form
$$\begin{pmatrix}
\begin{matrix}
1 & 0 & \cdots & 0 \\
0 & 1 & \cdots & 0 \\
\vdots & \vdots & \ddots & \vdots \\
0 & 0 & \cdots & 1 \\
\end{matrix}
& \cdots & * & \cdots & *\\
\vdots & \ddots & \vdots & \vdots & \vdots \\
p\mathbb{Z}_p & \cdots & 
\begin{matrix}
1 & 0 & \cdots & 0 \\
0 & 1 & \cdots & 0 \\
\vdots & \vdots & \ddots & \vdots \\
0 & 0 & \cdots & 1 \\
\end{matrix}
& \cdots & * \\
\end{pmatrix}.$$
Similarly, we can eliminate those $p\mathbb{Z}_p$ entries and obtain a basis of $\mathcal{L}(b_1,\dots,b_{s_{l-1}})$ of the form
$$\begin{pmatrix}
\begin{matrix}
1 & 0 & \cdots & 0 \\
0 & 1 & \cdots & 0 \\
\vdots & \vdots & \ddots & \vdots \\
0 & 0 & \cdots & 1 \\
\end{matrix}
& \cdots & * & \cdots & *\\
\vdots & \ddots & \vdots & \vdots & \vdots \\
0 & \cdots & 
\begin{matrix}
1 & 0 & \cdots & 0 \\
0 & 1 & \cdots & 0 \\
\vdots & \vdots & \ddots & \vdots \\
0 & 0 & \cdots & 1 \\
\end{matrix}
& \cdots & * \\
\end{pmatrix}.$$\par
When $t=l$, if there are coefficients $\overline{b_{s_{l-1}+1}},\overline{b_{s_{l-1}+2}},\dots,\overline{b_{s_l}}\in\mathbb{Z}_p/p\mathbb{Z}_p$, not all $\overline{0}$, such that
$$\overline{b_{s_{l-1}+1}}u_{s_{l-1}+1}+\overline{b_{s_{l-1}+2}}u_{s_{l-1}+2}+\cdots+\overline{b_{s_l}}u_{s_l}=\overline{0},$$
then the coefficients of $\alpha_j$ $({s_{l-1}+1}\le j\le {s_l})$ in
$$v=\sum^{s_l}_{i={s_{l-1}+1}}{b_i\beta_i}$$
all belong to $p\mathbb{Z}_p$. Combining with the previous conclusion, we obtain that the coefficients of $\alpha_j$ $(1\le j\le {s_l})$ in the above sum all belong to $p\mathbb{Z}_p$. Then we can use $b_1,\dots,b_{s_{l-1}}$ to eliminate the first $s_{l-1}$ columns of $v$. Meanwhile, the $(s_{l-1}+1)$th to the $s_l$th columns of $v$ still belong to $p\mathbb{Z}_p$, i.e., the following matrix
$$\begin{pmatrix}
\begin{matrix}
1 & 0 & \cdots & 0 \\
0 & 1 & \cdots & 0 \\
\vdots & \vdots & \ddots & \vdots \\
0 & 0 & \cdots & 1 \\
\end{matrix}
& \cdots & * & * & \cdots & * & * & \cdots & * \\
\vdots & \ddots & \vdots & \vdots & \ddots & \vdots & \vdots & \ddots & \vdots \\
0 & \cdots & 
\begin{matrix}
1 & 0 & \cdots & 0 \\
0 & 1 & \cdots & 0 \\
\vdots & \vdots & \ddots & \vdots \\
0 & 0 & \cdots & 1 \\
\end{matrix}
& *  & \cdots & * & * & \cdots & * \\
p\mathbb{Z}_p & \cdots & p\mathbb{Z}_p & p\mathbb{Z}_p & \cdots & p\mathbb{Z}_p & * & \cdots & * \\
\end{pmatrix}$$
becomes
$$\begin{pmatrix}
\begin{matrix}
1 & 0 & \cdots & 0 \\
0 & 1 & \cdots & 0 \\
\vdots & \vdots & \ddots & \vdots \\
0 & 0 & \cdots & 1 \\
\end{matrix}
& \cdots & * & * & \cdots & * & * & \cdots & * \\
\vdots & \ddots & \vdots & \vdots & \ddots & \vdots & \vdots & \ddots & \vdots \\
0 & \cdots & 
\begin{matrix}
1 & 0 & \cdots & 0 \\
0 & 1 & \cdots & 0 \\
\vdots & \vdots & \ddots & \vdots \\
0 & 0 & \cdots & 1 \\
\end{matrix}
& *  & \cdots & * & * & \cdots & * \\
0 & \cdots & 0 & c_{s_{l-1}+1} & \cdots & c_{s_l} & c_{s_l+1} & \cdots & c_n \\
\end{pmatrix},$$
where $c_{s_{l-1}+1},\dots,c_{s_l}\in p\mathbb{Z}_p$. Hence there exist $b_1,\dots,b_{s_{l-1}}\in\mathbb{Z}_p$ such that
$$N\left(\sum^{s_{l-1}}_{i=1}{b_i\beta_i}+v\right)=N\left(\sum^{n}_{i=s_{l-1}+1}{c_i\alpha_i}\right)=\max_{s_{l-1}+1\le i\le n}{N(c_i\alpha_i)}<N(\alpha_{s_l})=N(\beta_{s_l}).$$
On the other hand, 
$$N\left(\sum^{s_{l-1}}_{i=1}{b_i\beta_i}+v\right)=N\left(\sum^{s_l}_{i=1}{b_i\beta_i}\right)=\max_{1\le i\le s_l}{N(b_i\beta_i)}\ge N(\beta_{s_l}).$$
This is a contradiction. The proof of the claim is complete.\par
Let us come back to the proof of the theorem. Notice that $N(\alpha_m)=N(\beta_m)\ge N(\beta_{m+1})>N(\alpha_{m+1})$, hence $s_k=m$. According to the above results, there exist coefficients $b_1,\dots,b_m\in\mathbb{Z}_p$ such that
$$N\left(\sum^{m}_{i=1}{b_i\beta_i}+\beta_{m+1}\right)=N\left(\sum^{n}_{i=m+1}{c_i\alpha_i}\right)=\max_{m+1\le i\le n}{N(c_i\alpha_i)}\le N(\alpha_{m+1})<N(\beta_{m+1}).$$
On the other hand,
$$N\left(\sum^{m}_{i=1}{b_i\beta_i}+\beta_{m+1}\right)=\max{\{N(b_1\beta_1),\dots,N(b_m\beta_m),N(\beta_{m+1})\}}\ge N(\beta_{m+1}).$$
This is a contradiction. Hence we must have $N(\alpha_{m+1})=N(\beta_{m+1})$. The proof of the theorem is complete.
\end{proof}

\subsection{Successive Maxima}

Now we can define the successive maxima.
\begin{definition}[successive maxima]
Let $V$ be a vector space over $\mathbb{Q}_p$ of finite dimension, and let $N$ be a norm on $V$. Let $\mathcal{L}$ be a $p$-adic lattice of rank $n$ in $V$. Let $\alpha_1,\dots,\alpha_n$ be an $N$-orthogonal basis of $\mathcal{L}$ such that $N(\alpha_{1})\ge\cdots\ge N(\alpha_{n})$. The $i$th successive maxima of $\mathcal{L}$ respect to norm $N$ is
$$\tilde{\lambda}_i(\mathcal{L}):=N(\alpha_i).$$
\end{definition}
By Main Theorem \ref{th-2.1}, it is well defined. The successive maxima can also be defined through an iterative process. First, we choose the longest vector $\alpha_1\in\mathcal{L}$ and define $\tilde{\lambda}_1(\mathcal{L})=N(\alpha_1)$. Suppose that $\tilde{\lambda}_1(\mathcal{L})=N(\alpha_1),\dots,\tilde{\lambda}_i(\mathcal{L})=N(\alpha_i)$ are defined. Then we choose the longest vector $\alpha_{i+1}\in\mathcal{L}$ such that $\alpha_{i+1}\mathbb{Z}_p$ is $N$-orthogonal to $\alpha_1\mathbb{Z}_p+\cdots+\alpha_i\mathbb{Z}_p$ and define $\tilde{\lambda}_{i+1}(\mathcal{L})=N(\alpha_{i+1})$.\par
We have the following property of the successive maxima.

\begin{proposition}\label{pr-2.3}
Let $V$ be a vector space over $\mathbb{Q}_p$ of finite dimension, and let $N$ be a norm on $V$. Let $\mathcal{L}=\mathcal{L}(\alpha_1,\dots,\alpha_n)$ be a $p$-adic lattice of rank $n$ in $V$ such that $N(\alpha_{1})\ge\cdots\ge N(\alpha_{n})$. Then $\tilde{\lambda}_1(\mathcal{L})=N(\alpha_1)$ and $\tilde{\lambda}_i(\mathcal{L})\le N(\alpha_i)$ for $2\le i\le n$.
\begin{proof}
Since the longest vector in $\mathcal{L}$ is $\alpha_1$, we have $N(\alpha_1)=\tilde{\lambda}_1(\mathcal{L})$ by definition. We can use the algorithm in Main Theorem \ref{th-1.5} to find an $N$-orthogonal basis $\beta_1,\dots,\beta_n$ of $\mathcal{L}$. Notice that the order of $\alpha_1,\dots,\alpha_n$ may be changed during the orthogonalization process. We can initially arrange their order in such a way that this sequence remains unchanged throughout the orthogonalization process, as we mentioned in the proof of Main Theorem \ref{th-1.5}. Then we have
$$N(\beta_{1})\ge\cdots\ge N(\beta_{n})$$
and
$$N(\beta_{i})\le N(\alpha_{i})$$
for $1\le i\le n$. Since
$$N(\alpha_{1})\ge\cdots\ge N(\alpha_{n}),$$
we have
$$N(\beta_{i})\le N(\alpha_{i})$$
for $1\le i\le n$. Finally, by definition, we have $N(\beta_{i})=\tilde{\lambda}_i(\mathcal{L})$ and hence $\tilde{\lambda}_i(\mathcal{L})\le N(\alpha_i)$ for $2\le i\le n$.
\end{proof}
\end{proposition}

This bound is tight, as demonstrated by the fact that if $\alpha_1,\dots,\alpha_n$ is an $N$-orthogonal basis of $\mathcal{L}$, then, by definition, $\tilde{\lambda}_i(\mathcal{L})=N(\alpha_i)$ for $1\le i\le n$. The subsequent corollary establishes the converse of this statement.

\begin{corollary}\label{co-2.4}
Let $V$ be a vector space over $\mathbb{Q}_p$ of finite dimension, and let $N$ be a norm on $V$. Let $\mathcal{L}=\mathcal{L}(\alpha_1,\dots,\alpha_n)$ be a $p$-adic lattice of rank $n$ in $V$ such that $\tilde{\lambda}_i(\mathcal{L})=N(\alpha_i)$ for $1\le i\le n$. Then $\alpha_1,\dots,\alpha_n$ is an $N$-orthogonal basis of $\mathcal{L}$.
\begin{proof} Suppose that $\alpha_1,\dots,\alpha_n$ is not an $N$-orthogonal basis of $\mathcal{L}$. Then there is a subscript $i^{\prime}$ such that $\alpha_1,\dots,\alpha_{i^{\prime}-1}$ is an $N$-orthogonal basis of $\mathcal{L}(\alpha_1,\dots,\alpha_{i^{\prime}-1})$ while $\alpha_1,\dots,\alpha_{i^{\prime}}$ is not an $N$-orthogonal basis of $\mathcal{L}(\alpha_1,\dots,\alpha_{i^{\prime}})$. Since $\tilde{\lambda}_i(\mathcal{L})=N(\alpha_i)$ for $1\le i\le n$, we have $N(\alpha_1)\ge\cdots\ge N(\alpha_n)$. Let $w_{0}\in\mathcal{L}(\alpha_{1},\dots,\alpha_{i^{\prime}-1})$ be such that
$$N(\alpha_{i^{\prime}}+w_{0})=\min\{N(\alpha_{i^{\prime}}+w):w\in\mathcal{L}(\alpha_{1},\dots,\alpha_{i^{\prime}-1})\}.$$
Then, we have $N(\alpha_{i^{\prime}}+w_{0})<N(\alpha_{i^{\prime}})$, otherwise we can take $w_{0}=0$ and $\alpha_1,\dots,\alpha_{i^{\prime}}$ is an $N$-orthogonal basis of $\mathcal{L}(\alpha_1,\dots,\alpha_{i^{\prime}})$ by Lemma \ref{le-1.4}. Replace $\alpha_{i^{\prime}}$ by $\alpha_{i^{\prime}}+w_{0}$ in the basis, we have
$$\mathcal{L}=\mathcal{L}(\alpha_1,\dots,\alpha_{i^{\prime}-1},\alpha_{i^{\prime}}+w_{0},\alpha_{i^{\prime}+1},\dots,\alpha_n).$$
Aussume that $N(\alpha_{j})\ge N(\alpha_{i^{\prime}}+w_{0})\ge N(\alpha_{j+1})$ for some $i^{\prime}\le j\le n$ (if $j=n$, then there is olny $N(\alpha_{n})\ge N(\alpha_{i^{\prime}}+w_{0})$). Then, by Proposition \ref{pr-2.3},
$$N(\alpha_{i^{\prime}+i})\ge\tilde{\lambda}_{i^{\prime}+i-1}(\mathcal{L})$$
for $1 \le i\le j-i^{\prime}$, and
$$N(\alpha_{i^{\prime}}+w_{0})\ge\tilde{\lambda}_{j}(\mathcal{L}).$$
On the other hand, since
$$N(\alpha_{i^{\prime}+i})=\tilde{\lambda}_{i^{\prime}+i}(\mathcal{L})\le\tilde{\lambda}_{i^{\prime}+i-1}(\mathcal{L}),$$
we have
$$\tilde{\lambda}_{i^{\prime}+i}(\mathcal{L})=\tilde{\lambda}_{i^{\prime}+i-1}(\mathcal{L})$$
for $1 \le i\le j-i^{\prime}$. Therefore,
$$N(\alpha_{i^{\prime}}+w_{0})\ge\tilde{\lambda}_{j}(\mathcal{L})=\tilde{\lambda}_{i^{\prime}}(\mathcal{L})=N(\alpha_{i^{\prime}})>N(\alpha_{i^{\prime}}+w_{0}),$$
which is a contradiction.
\end{proof}
\end{corollary}

If we know the successive maxima of a $p$-adic lattice $\mathcal{L}$, then we can compute all possible norms of vectors in $\mathcal{L}$, which are
$$N(\mathcal{L})=\{N(v):v\in\mathcal{L}\}=\{p^{-i}\tilde{\lambda}_j(\mathcal{L}):i\in\mathbb{Z}_{\ge0},j=1,\dots,n\}\cup\{0\}.$$
Hence we can use successive maxima to solve the LVP in $p$-adic lattices.\par

%

\subsection{Escape Distance}

In the context of Euclidean lattices, the concept of covering radius represents the maximal distance between a full-rank lattice and a point lying outside it. By substituting ``maximal'' with ``minimal'', we obtain the $p$-adic analogue of this definition.

\begin{definition}[escape distance]
For a full rank $p$-adic lattice $\mathcal{L}$, define the escape distance of $\mathcal{L}$ as
$$\mu(\mathcal{L}):=\min_{x\in V\setminus\mathcal{L}}{{\rm dist}(x,\mathcal{L})}.$$
\end{definition}

Here is a toy example.

\begin{example}
Let $V=\mathbb{Q}_p$ and $\mathcal{L}=\mathbb{Z}_p$. Then $\mu(\mathcal{L})=N(p^{-1})=p$.
\end{example}

We can determine the escape distance by the successive maxima.

\begin{theorem}
For a full rank $p$-adic lattice $\mathcal{L}$ of rank $n$, we have $\mu(\mathcal{L})=p\tilde{\lambda}_n(\mathcal{L})$.
\begin{proof}
Let $\alpha_1,\dots,\alpha_n$ be an $N$--orthogonal basis of $\mathcal{L}$ such that $N(\alpha_i)=\tilde{\lambda}_i(\mathcal{L})$ for $i=1,\dots,n$. Then for any lattice vector $v\in\mathcal{L}$ and vector $x\in V\setminus\mathcal{L}$, we can write $v=\sum_{i=1}^{n}{a_i\alpha_i}$ and $x=\sum_{i=1}^{n}{b_i\alpha_i}$ where $a_i\in\mathbb{Z}_p$ and $b_i\in\mathbb{Q}_p$ for $i=1,\dots,n$. Since $x\not\in\mathcal{L}$, at least one $b_i$ is not in $\mathbb{Z}_p$. Suppose that $b_{i^{\prime}}\not\in\mathbb{Z}_p$. Then we have
\begin{equation*}
\begin{split}
N(x-v)&=\max_{1\le i\le n}{N\big((b_i-a_i)\alpha_i\big)}\ge N\big((b_{i^{\prime}}-a_{i^{\prime}})\alpha_{i^{\prime}}\big)\\
&=N(b_{i^{\prime}}\alpha_{i^{\prime}})=\left|b_{i^{\prime}}\right|_p\cdot N(\alpha_{i^{\prime}})\\
&\ge pN(\alpha_n).
\end{split}
\end{equation*}
The equality holds when $x=p^{-1}\alpha_n$. Hence $\mu(\mathcal{L})=pN(\alpha_n)=p\tilde{\lambda}_n(\mathcal{L})$.
\end{proof}
\end{theorem}

\section{Transformation of $N$-orthogonal Basis}\label{se-3}

The following theorem characterizes transformations between $N$-orthogonal bases of a $p$-adic lattice.

\begin{theorem}
Let $V$ be a vector space over $\mathbb{Q}_p$ of finite dimension, and let $N$ be a norm on $V$. Let $\mathcal{L}=\mathcal{L}(\alpha_1,\dots,\alpha_n)$ be a $p$-adic lattice of rank $n$ in $V$ with an $N$-orthogonal basis $\alpha_1,\dots,\alpha_n$. Then, $\beta_1,\dots,\beta_n$ is also an $N$-orthogonal basis of $\mathcal{L}$ if and only if it can be obtained from $\alpha_1,\dots,\alpha_n$ by the following operations:
\begin{enumerate}
\item $\alpha_i\leftarrow k\alpha_i$ for some $k\in\mathbb{Z}_p\setminus p\mathbb{Z}_p$,
\item $\alpha_i\leftrightarrow\alpha_j$,
\item $\alpha_i\leftarrow\alpha_i+k\alpha_j$ for some $k\in\mathbb{Z}_p$ such that $N(k\alpha_j)\le N(\alpha_i)$,
\end{enumerate}
\begin{proof}
First, we prove the necessity. Clearly, $\beta_1,\dots,\beta_n$ is still a basis of $\mathcal{L}$. Since
$$N(\alpha_i+k\alpha_j)=\max{\{N(\alpha_i),N(k\alpha_j)\}}=N(\alpha_i),$$
the sorted norm sequence keeps the same after one step of these operations. By Corollary \ref{co-2.4}, after one step of these operations, the new basis is an $N$-orthogonal basis. Therefore, we conclude by induction that $\beta_1,\dots,\beta_n$ is an $N$-orthogonal basis of $\mathcal{L}$.\par
Now, we prove the sufficiency. We can write
$$\beta_i=\sum^{n}_{j=1}{a_{ij}\alpha_j},$$
where $a_{ij}\in\mathbb{Z}_p$ for $1\le i,j\le n$. Let the matrix $A=(a_{ij})$ be the coefficient matrix. Our aim is to demonstrate that matrix $A$ can be converted to an identity matrix through these three operations. The transformation is consistent with the steps employed in the proof of Main Theorem \ref{th-2.1}. Now we prove that each step in this process corresponds to one of the three defined operations.\par
First, when we try to obtain the small identity matrix in $(\mathbb{Z}_p/p\mathbb{Z}_p)^{*\times *}$, we use the following operations:
\begin{enumerate}
\item $\beta_i\leftarrow k\beta_i$ for some $k\in\mathbb{Z}_p\setminus p\mathbb{Z}_p$,
\item $\beta_i\leftrightarrow\beta_j$,
\item $\beta_i\leftarrow\beta_i+k\beta_j$ for some $k\in\mathbb{Z}_p$.
\end{enumerate}
Since $N(\beta_{s_{t-1}+1})=\cdots=N(\beta_{s_t})$, we have $N(k\beta_j)\le N(\beta_i)$ for all $k\in\mathbb{Z}_p$ and $s_{t-1}+1\le i,j\le s_t$. Therefore, these operations are exactly the same as those in this theorem.\par
Next, when we attempt to obtain the small identity matrix in $(\mathbb{Z}_p)^{*\times *}$ from the small identity matrix in $(\mathbb{Z}_p/p\mathbb{Z}_p)^{*\times *}$, we use the same kinds of operations as above.\par
Following this, we seek to eliminate $a_{ij}\in p\mathbb{Z}_p$ for $s_{t-1}+1\le i\le s_t$ and $1\le j\le s_{t-1}$. We need to prove that $N(a_{ij}\beta_j)\le N(\beta_i)$. If $N(a_{i^{\prime}j^{\prime}}\beta_{j^{\prime}})>N(\beta_{i^{\prime}})$, then $N(a_{i^{\prime}j^{\prime}}\alpha_{j^{\prime}})>N(\alpha_{i^{\prime}})$ and
$$N(\alpha_{i^{\prime}})=N(\beta_{i^{\prime}})=N\left(\sum^{n}_{j=1}{a_{i^{\prime}j}}\alpha_{j}\right)=\max_{1\le j\le n}{N(a_{i^{\prime}j}\alpha_{j})}\ge N(a_{i^{\prime}j^{\prime}}\alpha_{j^{\prime}})>N(\alpha_{i^{\prime}}).$$
This is a contradiction.\par
Now we obtain an upper triangle matrix with diagonals all equal to $1$. Since $N(k\beta_j)\le N(\beta_i)$ for all $k\in\mathbb{Z}_p$ and $1\le i\le j\le n$, we can use operation $3$ to transform it to an identity matrix.
\end{proof}
\end{theorem}

%

\section{Orthogonalization with Orthogonal Bases of the Vector Space}\label{se-4}

Efficiently solving the CVP in $p$-adic lattices remains an open challenge without the aid of an $N$-orthogonal basis for the vector space, and the algorithm presented in the proof of Main Theorem \ref{th-1.5} is not considered efficient. However, possessing an $N$-orthogonal basis of the vector space enables us to efficiently determine an $N$-orthogonal basis of a $p$-adic lattice in this vector space.\par
Let $V$ be a vector space over $\mathbb{Q}_p$ of finite dimension $n>0$, and let $N$ be a norm on $V$. Let $e_1,\dots,e_n$ be an $N$-orthogonal basis of $V$ over $\mathbb{Q}_p$. Let $\mathcal{L}=\mathcal{L}(\alpha_1,\dots,\alpha_m)$ be a $p$-adic lattice of rank $m\le n$ in $V$. Then we can write
$$\alpha_i=\sum^{n}_{j=1}{a_{ij}e_j},$$
where $a_{ij}\in\mathbb{Q}_p$ for $i=1,\dots,m$ and $j=1,\dots,n$. We have the following orthogonalization process.\par\vskip 10pt

{\bf Algorithm} (orthogonalization with orthogonal bases of the vector space).\par
{\bf Input:} an $N$-orthogonal basis $e_1,\dots,e_n$ of $V$, a $p$-adic lattice $\mathcal{L}=\mathcal{L}(\alpha_1,\dots,\alpha_m)$ in $V$.\par
{\bf Output:} an $N$-orthogonal basis of $\mathcal{L}$.
\begin{enumerate}
\item for $i=1$ to $m$ do:
\item\quad rearrange $\alpha_i,\dots,\alpha_m$ such that  $N(\alpha_i)=\max_{i\le k\le m}{N(\alpha_k)}$,
\item\quad rearrange $e_i,\dots,e_n$ such that  $N(a_{ii}e_i)=\max_{i\le j\le m}{N(a_{ij}e_j)}$,
\item\quad for $l=i+1$ to $m$ do:
\item\quad\quad $\alpha_l\leftarrow\alpha_l-\frac{a_{li}}{a_{ii}}\alpha_i.$
\end{enumerate}\par
Return $(\alpha_1,\dots,\alpha_m)$.\par\vskip 10pt

This algorithm runs in polynomial time in the input size if we can compute efficiently the norm $N(v)$ of any vector $v\in V$. Let us now highlight some crucial observations regarding this algorithm. In the first iteration of the outer loop, we have $N(\alpha_1)=\max_{1\le k\le m}{N(\alpha_k)}$ and $N(a_{11}e_1)=\max_{1\le j\le m}{N(a_{1j}e_j)}$. Then we use $a_{11}$ to eliminate $a_{l1}$ for $2\le l\le m$. After the first iteration of the outer loop, the matrix of $\alpha_1,\dots,\alpha_m$ is of the form
$$\begin{pmatrix}
a_{11} & a_{12} & \cdots & a_{1n} \\
0 & a^{\prime}_{22} & \cdots & a^{\prime}_{2n} \\
\vdots & \vdots & \ddots & \vdots \\
0 & a^{\prime}_{m2} & \cdots & a^{\prime}_{mn} \\
\end{pmatrix}.$$
Moreover, $\alpha_1$ and $e_1$ are fixed during the rest of the algorithm. Similarly, in the $i$th iteration of the outer loop, the entries in the $(i+1)$th to the $m$th row of the $i$th column are eliminated. Also, $\alpha_i$ and $e_i$ are fixed after the $i$th iteration of the outer loop. Finally, we obtain a matrix of the form
$$\begin{pmatrix}
b_{11} & b_{12} & \cdots & b_{1m} & \cdots & b_{1n}\\
0 & b_{22} & \cdots & b_{2m} & \cdots & b_{2n}\\
\vdots & \vdots & \ddots & \vdots & \ddots & \vdots\\
0 & 0 & \cdots & b_{mm} & \cdots & b_{mn}\\
\end{pmatrix}.$$
Denote the row vectors of this matrix by $\beta_i$, $1\le i\le m$. Let us prove that it is an $N$-orthogonal basis of $\mathcal{L}$.
\begin{theorem}
The algorithm above outputs an $N$-orthogonal basis of $\mathcal{L}$.
\begin{proof}
We first prove that the algorithm outputs a basis of $\mathcal{L}$. In the $i$th iteration of the outer loop, we have $N(\alpha_i)=\max_{i\le k\le m}{N(\alpha_k)}$, $N(a_{ii}e_i)=\max_{i\le j\le m}{N(a_{ij}e_j)}$ and $a_{ij}=0$ for $1\le j\le i-1$. Since $N(\alpha_i)=\max_{1\le j\le m}{N(a_{ij}e_j)}$, we can imply that
$$N(a_{ii}e_i)=N(\alpha_i)\ge N(\alpha_l)\ge N(a_{li}e_i)$$
for $i+1\le l\le m$. Hence $N(a_{ii})\ge N(a_{li})$ and $\frac{a_{li}}{a_{ii}}\in\mathbb{Z}_p$. Therefore, the algorithm outputs a basis of $\mathcal{L}$.\par
Now we prove by induction that, after the $i$th iteration of the outer loop, $\beta_1,\dots,\beta_i$ constitutes an $N$-orthogonal basis of $\mathcal{L}(\beta_1,\dots,\beta_i)$. When $i=1$ the conclusion is trivial. Suppose that the conclusion holds for $i=t$. Then $\beta_1,\dots,\beta_t$ is an $N$-orthogonal basis of $\mathcal{L}(\beta_1,\dots,\beta_t)$ by induction hypothesis. Furthermore, we observe $N(\beta_1)\ge\cdots\ge N(\beta_t)\ge N(\beta_{t+1})$, as we consistently select the longest vector at the beginning of each outer loop, and the operations within the inner loop can not  augment the length of $\alpha_l$ beyond that of $\alpha_i$. Therefore, according to Lemma \ref{le-1.4}, we just need to prove that
$$N(\beta_{t+1})=\min\left\{N(\beta_{t+1}+w):w\in\mathcal{L}(\beta_1,\dots,\beta_t)\right\}.$$
If there is
$$w=\sum^{t}_{i=1}{c_i\beta_i}$$
where $c_i\in\mathbb{Z}_p$, $1\le i\le t$, such that $N(\beta_{t+1}+w)<N(\beta_{t+1})$, then we must have $N(w)=N(\beta_{t+1})$. The $e_j$ coordinate of $\beta_{t+1}+w$, denoted by $d_j$, is
$$\sum^{t}_{i=1}{c_ib_{ij}}=\sum^{j}_{i=1}{c_ib_{ij}}$$
for $1\le j\le t$ and
$$\sum^{t}_{i=1}{c_ib_{ij}}+b_{(t+1)j}$$
for $t+1\le j\le m$. Since
$$\max_{1\le j\le m}{N(d_je_j)}=N(\beta_{t+1}+w)<N(\beta_{t+1})=\max_{t+1\le j\le m}{N\left(b_{(t+1)j}e_j\right)},$$
then for each $t+1\le j_0\le m$ such that 
$$N\left(b_{(t+1)j_0}e_{j_0}\right)=\max_{t+1\le j\le m}{N\left(b_{(t+1)j}e_j\right)},$$
we must have
$$N(d_{j_0}e_{j_0})=N\left(\sum^{t}_{i=1}{c_ib_{ij_0}}e_{j_0}+b_{(t+1)j_0}e_{j_0}\right)<N\left(b_{(t+1)j_0}e_{j_0}\right).$$
Hence
$$N(\beta_{t+1})=N\left(b_{(t+1)j_0}e_{j_0}\right)=N\left(\sum^{t}_{i=1}{c_ib_{ij_0}}e_{j_0}\right)\le\max_{1\le i\le t}{N(c_ib_{ij_0}e_{j_0})}.$$
Suppose that $i_0$ is the first subscript such that $N(c_{i_0}b_{i_0j_0}e_{j_0})\ge N(\beta_{t+1})$. Then for $1\le i<i_0$, we have
$$N(c_ib_{ij_0}e_{j_0})<N(\beta_{t+1})\le N(c_{i_0}b_{i_0j_0}e_{j_0}).$$
Hence $N(c_ib_{ii_0}e_{i_0})<N(c_{i_0}b_{i_0i_0}e_{i_0})$ for $1\le i<i_0$. Therefore, the norm of the $e_{i_0}$ coordinate of $\beta_{t+1}+w$ is
$$N(d_{i_0}e_{i_0})=N\left(\sum^{i_0}_{i=1}{c_ib_{ii_0}e_{i_0}}\right)=N(c_{i_0}b_{i_0i_0}e_{i_0})\ge N(\beta_{t+1}).$$
This is a contradiction.
\end{proof}
\end{theorem}

Here is a toy example to explain the above algorithm.

\begin{example}
Keep the notation in Example \ref{eg-1.6}. Here, we do not solve CVP instances. Instead, we use the $N$-orthogonal basis $1,\zeta,\zeta^2,\zeta^3$.\par
Since $N(1)=N(1+2\zeta)=1$ are the longest, this time we can choose $1+2\zeta$ first. Since $N(1)>N(2\zeta)$, we eliminate the $1$ coordinate of the vectors $1$ and $2+8\zeta+16\zeta^2+16\zeta^3$ by the vector $1+2\zeta$. We obtain the vectors $-2\zeta$ and $4\zeta+16\zeta^2+16\zeta^3$, respectively.\par
Next, since $N(-2\zeta)>N(4\zeta+16\zeta^2+16\zeta^3)$, we eliminate the $\zeta$ coordinate of the vector $4\zeta+16\zeta^2+16\zeta^3$ by the vector $-2\zeta$. We obtain the vector $16\zeta^2+16\zeta^3$. Hence $1+2\zeta,-2\zeta,16\zeta^2+16\zeta^3$ is an $N$-orthogonal basis of $\mathcal{L}$.
\end{example}

\section{Solving the CVP and the LVP with Orthogonal Bases}\label{se-5}

Given our ability to perform the orthogonalization process using an $N$-orthogonal basis of the vector space, a pertinent question arises: Can we solve the CVP utilizing an $N$-orthogonal basis of the vector space? The solution presented in Theorem 3.6 of \cite{ref-2} addresses the CVP specifically with the help of orthogonal bases. However, it requires the lattice basis to be included in the $N$-orthogonal basis of the vector space, a more restrictive condition than our current assumption.\par
In this section, we introduce an algorithm for solving the CVP with the help of an $N$-orthogonal basis of the vector space, and we assert that the orthogonalization process and the CVP are polynomially equivalent.\par
Let $V$ be a vector space over $\mathbb{Q}_p$ of finite dimension $n>0$, and let $N$ be a norm on $V$. Let $e_1,\dots,e_n$ be an $N$-orthogonal basis of $V$ over $\mathbb{Q}_p$. Let $\mathcal{L}=\mathcal{L}(\alpha_1,\dots,\alpha_m)$ be a $p$-adic lattice of rank $m\le n$ in $V$. Then we can write
$$\alpha_i=\sum^{n}_{j=1}{a_{ij}e_j},$$
where $a_{ij}\in\mathbb{Q}_p$ for $i=1,\dots,m$ and $j=1,\dots,n$. Let $t\in V$ be a target vector. We present the following algorithm to solve the CVP with orthogonal bases.\par\vskip 10pt

{\bf Algorithm} (CVP with orthogonal bases).\par
{\bf Input:} an $N$-orthogonal basis $e_1,\dots,e_n$ of $V$, a $p$-adic lattice $\mathcal{L}=\mathcal{L}(\alpha_1,\dots,\alpha_m)$ in $V$, a target vector $t\in V$.\par
{\bf Output:} a closest lattice vector $v$ of $t$.
\begin{enumerate}
\item $v\leftarrow0$, $s\leftarrow t$, write $s=\sum^{n}_{j=1}{s_je_j}$,
\item for $i=1$ to $m$ do:
\item\quad rearrange $\alpha_i,\dots,\alpha_m$ such that  $N(\alpha_i)=\max_{i\le k\le m}{N(\alpha_k)}$,
\item\quad if $N(s)>N(\alpha_i)$ then break,
\item\quad rearrange $e_i,\dots,e_n$ such that  $N(a_{ii}e_i)=\max_{i\le j\le m}{N(a_{ij}e_j)}$,
\item\quad $s\leftarrow s-\frac{s_i}{a_{ii}}\alpha_i$, $v\leftarrow v+\frac{s_i}{a_{ii}}\alpha_i$,
\item\quad if $s=0$ then break,
\item\quad for $l=i+1$ to $m$ do:
\item\quad\quad $\alpha_l\leftarrow\alpha_l-\frac{a_{li}}{a_{ii}}\alpha_i$.
\end{enumerate}\par
Return $v$.\par\vskip 10pt

This algorithm runs in polynomial time in the input size if we can compute efficiently the norm $N(v)$ of any vector $v\in V$. Let us prove that $v$ is a closest lattice vector of $t$.

\begin{theorem}\label{th-5.1}
The algorithm above outputs a closest vector of $t$.
\begin{proof}
We first prove that $v\in\mathcal{L}$. If the algorithm performs step $6$ in the $i$th iteration of the outer loop, then we must have $N(s_ie_i)\le N(s)\le N(\alpha_i)=N(a_{ii}e_i)$. Hence $\frac{s_i}{a_{ii}}\in\mathbb{Z}_p$ and $v$ is a lattice vector.\par
We now prove that $v$ is a closest lattice vector of $t$. If the algorithm breaks at step $7$, i.e., $s=0$, then we have $t\in\mathcal{L}$ and $v=t$ is the closest lattice vector. Now assume $t\not\in\mathcal{L}$, then the algorithm never breaks at step $7$. Suppose that the algorithm breaks in the $i_0$th iteration of the outer loop (if it never breaks, then set $i_0=m+1$). Notice that steps $2,3,5,8,9$ are exactly steps $1$ to $5$ in the orthogonalization process in section \ref{se-4}. Moreover, this algorithm does the same thing to $s$ if $N(s)\le N(\alpha_i)$. Hence after the $i_0$th iteration of the outer loop, we have
$$N(t-v)=N(s)=\min\left\{N(s+w):w\in\mathcal{L}(\alpha_1,\dots,\alpha_{i_0-1})\right\}.$$
Since $N(s+w)\ge N(s)>N(\alpha_i)$ for all $w\in\mathcal{L}(\alpha_1,\dots,\alpha_{i_0-1})$ and $i_0\le i\le m$, we have $N(s+w+w^{\prime})= N(s+w)\ge N(s)$ for all $w\in\mathcal{L}(\alpha_1,\dots,\alpha_{i_0-1})$ and $w^{\prime}\in\mathcal{L}(\alpha_{i_0},\dots,\alpha_m)$. Therefore, $v$ is a closest lattice vector of $t$.
\end{proof}
\end{theorem}

Finally, we can easily conclude that the orthogonalization process and the CVP are polynomially equivalent.

\begin{theorem}\label{th-5.2}
Finding an $N$-orthogonal basis of a $p$-adic lattice and solving the CVP in $p$-adic lattice are polynomially equivalent if we can compute efficiently the norm $N(v)$ of any vector $v\in V$.
\begin{proof}
If we can solve CVP-instances, then we can use the algorithm in the proof of Main Theorem \ref{th-1.5} to find an $N$-orthogonal basis of a $p$-adic lattice in polynomial time. Conversely, if we are able to find an $N$-orthogonal basis of a $p$-adic lattice, then we can choose a basis of the vector space and view it as a $p$-adic lattice. The $N$-orthogonal basis of this lattice is also an $N$-orthogonal basis of the vector space. Subsequently, according to Theorem \ref{th-5.1}, we can solve CVP-instances in polynomial time.
\end{proof}
\end{theorem}


Certainly, we can first compute an $N$-orthogonal basis of a $p$-adic lattice by the algorithm in Section \ref{se-4} and then use Theorem 3.5 of \cite{ref-2} to solve the LVP. Yet, if the goal is solely to identify a (second) longest vector, there is no need to compute the entire $N$-orthogonal basis of a $p$-adic lattice. Consequently, the algorithm can be simplified for this specific purpose.\par\vskip 10pt

{\bf Algorithm} (LVP with orthogonal bases).\par
{\bf Input:} an $N$-orthogonal basis $e_1,\dots,e_n$ of $V$, a $p$-adic lattice $\mathcal{L}=\mathcal{L}(\alpha_1,\dots,\alpha_m)$ in $V$.\par
{\bf Output:} a (second) longest vector of $\mathcal{L}$.
\begin{enumerate}
\item for $i=1$ to $m$ do:
\item\quad rearrange $\alpha_i,\dots,\alpha_m$ such that  $N(\alpha_i)=\max_{i\le k\le m}{N(\alpha_k)}$,
\item\quad if $i>1$ and $N(\alpha_{i-1})>N(\alpha_i)$ then break,
\item\quad rearrange $e_i,\dots,e_n$ such that  $N(a_{ii}e_i)=\max_{i\le j\le m}{N(a_{ij}e_j)}$,
\item\quad for $l=i+1$ to $m$ do:
\item\quad\quad $\alpha_l\leftarrow\alpha_l-\frac{a_{li}}{a_{ii}}\alpha_i$,
\item if $N(p\alpha_1)>N(\alpha_i)$ then $v\leftarrow p\alpha_1$,
\item else $v\leftarrow\alpha_i$.
\end{enumerate}\par
Return $v$.\par\vskip 10pt

This algorithm runs in polynomial time in the input size if we can compute efficiently the norm $N(v)$ of any vector $v\in V$. Let us prove that $v$ is a (second) longest vector of $\mathcal{L}$.
\begin{theorem}
The algorithm above outputs a (second) longest vector of $\mathcal{L}$.
\begin{proof}
Since the step 6 can not make $N(\alpha_l)$ greater than $N(\alpha_i)$, when the iteration of the outer loop ends or breaks, we have $N(\alpha_1)=\cdots=N(\alpha_{i-1})>N(\alpha_i)$ and $i>1$. Moreover, the norms of the remaining vectors in this $N$-orthogonal basis are all less than or equal to $N(\alpha_i)$. Hence $\tilde{\lambda}_1(\mathcal{L})=\cdots=\tilde{\lambda}_{i-1}(\mathcal{L})=N(\alpha_1)$ and $\tilde{\lambda}_i(\mathcal{L})=N(\alpha_i)$. Therefore, the potential candidates for the (second) longest vector have norm either $N(p\alpha_1)$ or $N(\alpha_i)$.
\end{proof}
\end{theorem}

\section{Conclusion}

$p$-adic lattices exhibit the distinctive property of possessing orthogonal bases, whereas Euclidean lattices, in general, lack such bases. This distinct characteristic of $p$-adic lattices may find applications in cryptography and various other areas.\par
The algorithms proposed in this paper rely on the orthogonal bases. It would be valuable to explore efficient algorithms for solving the LVP and the CVP without the help of orthogonal bases. Conversely, determining whether the LVP and the CVP in $p$-adic lattices are NP-complete could also be an interesting area for research.\par
There are many other important concepts in the realm of Euclidean lattices.  However, their $p$-adic analogues have yet to be defined. For instance, is there a $p$-adic analogue of the dual lattice in Euclidean lattices? Notice that the $p$-adic norm can not induce an inner product. The definition of the dual lattice in $p$-adic lattices may differ.  Additionally, are there any transference theorems in $p$-adic lattices to bound the successive maxima, akin to Minkowski's bound on the successive minima in Euclidean lattices? It is our belief that there is much work to be done in exploring these problems.

\section*{Acknowledgements}

This work was supported by National Natural Science Foundation of China(No. 12271517) and National Key R\&D Program of China(No. 2020YFA0712300).

\end{document}